\documentclass[]{iac}

\newtagform{brackets}{(}{)}
\usetagform{brackets}

\usepackage{graphicx}
\usepackage{amsmath}
\usepackage{bm}
\usepackage[version=4]{mhchem}
\usepackage[per-mode = symbol]{siunitx}  
\DeclareSIUnit\permille{\text{\textperthousand}}

\usepackage[mode=buildnew]{standalone}
\usepackage{tikz}

\usepackage{hyperref} 
\hypersetup{
    colorlinks = true,  
	allcolors = black  
}

\RequirePackage{color}
\usepackage{xcolor}
\usepackage{soul}

\usepackage[super]{nth}
\usepackage{mathtools}

\usepackage{enumitem}

\usepackage{footnote}
\usepackage{threeparttable, tablefootnote}
\usepackage{multirow}
\usepackage{booktabs}

\usepackage{caption}
\usepackage{subcaption}

\allowdisplaybreaks

\newcommand{\norm}[1]{\left\lVert#1\right\rVert}

\def \StageThreeOne {1}
\def \StageThreeTwo {2}
\def \StageThreeThree {3}
\def \CoastingThree {4}
\def \StageFourOne {5}
\def \CoastingFour {6}
\def \StageFourTwo {7}
\def \Return {8}
\def \PReturn {9}

\def \tf {t_f}

\begin{document}

\IACpaperyear{21}
\IACpapernumber{D2.3.8}
\IACconference{72}
\IAClocation{Dubai, United Arab Emirates}
\IACdate{25-29 October 2021}
\IACcopyrightB{2021}{Mr. Boris Benedikter}

\title{Autonomous Upper Stage Guidance with Robust Splash-Down Constraint}

\IACauthor{Boris~Benedikter$^{\text{a}\ast}$, }{$^\text{a}$\textit{PhD Student, Department of Mechanical and Aerospace Engineering, Sapienza University of Rome, Via Eudossiana 18, 00184, Rome, Italy}, \uline{boris.benedikter@uniroma1.it}}

\IACauthor{Alessandro~Zavoli$^{\text{b}}$, }{$^\text{b}$\textit{Research Assistant, Department of Mechanical and Aerospace Engineering, Sapienza University of Rome, Via Eudossiana 18, 00184, Rome, Italy}, \uline{alessandro.zavoli@uniroma1.it}}

\IACauthor{Guido~Colasurdo$^{\text{c}}$, }{$^\text{c}$\textit{Full Professor, Department of Mechanical and Aerospace Engineering, Sapienza University of Rome, Via Eudossiana 18, 00184, Rome, Italy}, \uline{guido.colasurdo@uniroma1.it}}

\IACauthor{Simone~Pizzurro$^{\text{d}}$, }{$^\text{d}$\textit{Research Fellow, Launchers and Space Transportation Department, Italian Space Agency, Via del Politecnico snc, 00133, Rome, Italy}, \uline{simone.pizzurro@est.asi.it}}

\IACauthor{Enrico~Cavallini$^{\text{e}}$}{$^\text{e}$\textit{Head of Space Transportation Programs Office, Space Transportation, Space Infrastructures, and In-Orbit Servicing Department, Italian Space Agency, Via del Politecnico snc, 00133, Rome, Italy}, \uline{enrico.cavallini@asi.it}}

\abstract{%
This paper presents a novel algorithm, based on model predictive control (MPC), for the optimal guidance of a launch vehicle upper stage.
The proposed strategy not only maximizes the performance of the vehicle and its robustness to external disturbances, but also robustly enforces the splash-down constraint.
Indeed, uncertainty on the engine performance, and in particular on the burn time, could lead to a large footprint of possible impact points, which may pose a concern if the reentry points are close to inhabited regions. 
Thus, the proposed guidance strategy incorporates a neutral axis maneuver (NAM) that minimizes the sensitivity of the impact point to uncertain engine performance.
Unlike traditional methods to design a NAM, which are particularly burdensome and require long validation and verification tasks, the presented MPC algorithm autonomously determines the neutral axis direction by repeatedly solving an optimal control problem (OCP) with two return phases, a nominal and a perturbed one, constrained to the same splash-down point.
The OCP is transcribed as a sequence of convex problems that quickly converges to the optimal solution, thus allowing for 
high MPC update frequencies.
Numerical results assess the robustness and performance of the proposed algorithm via extensive Monte Carlo campaigns.
}

\maketitle

\section{Introduction}
Launch vehicle dynamics are subject to significant uncertainties, due to hard-to-predict aerodynamic coefficients, scattering of propulsion system performance, and sudden variations of the local environment, to name a few.
Therefore, a robust guidance algorithm that ensures the accurate injection of the payload into the desired orbit while meeting the numerous system requirements, even in presence of model uncertainties, is an absolute necessity for the success of a launch vehicle (LV) ascent mission.
Also, since launch vehicles, as most aerospace systems, must operate at the boundaries of their performance envelope to minimize the overall mission cost, their guidance should be \emph{optimal}, meaning that the vehicle should fly trajectories that are not only safe but also convenient (e.g., in terms of propellant consumption, to maximize the carrying capacity and the responsiveness in case of off-nominal conditions).

Model predictive control (MPC) is one of the few control synthesis strategies that can optimize the system performance while systematically accomodating mission constraints \cite{kouvaritakis2016model}.
Specifically, MPC consists of solving repeatedly an optimal control problem (OCP), updated with the onboard system measurements, and implementing the computed optimal control law in the time frames between the optimization procedures.
The closed-loop architecture inherently provides robustness to model uncertainties and in-flight disturbances, as the continuous update of the optimal path compensates for deviations from the nominal one.

Among the aerospace community, the interest toward 
MPC has grown constantly over the years \cite{eren2017model}, featuring successful employment over diverse problems, including rocket landing \cite{pascucci2015model,wang2019optimal}, spacecraft landing \cite{carson2006model, reynolds2017small}, rendezvous \cite{hartley2015tutorial,weiss2015model}, and formation flying \cite{keviczky2008decentralized, morgan2014model}, to name a few.
Indeed, since aerospace vehicles are typically subject to tight mission requirements, equipped with limited computational resources, and often operate in uncertain conditions, MPC represents an uniquely effective tool to control such systems.
The earliest applications of MPC were limited to relatively simple problems, as the computational cost associated with the solution of the OCP can impair the effectiveness of the MPC scheme, which can compensate for random disturbances only if the update frequency is sufficiently high.
For this reason, applying the appealing properties of MPC to nonlinear systems has always been a challenging task.
However, in the last decades, the advances in computing hardware and optimization theory gave a new drive to \emph{nonlinear MPC} \cite{mayne2000nonlinear}, which still poses serious challenges in terms of robustness guarantees and computational efficiency, but ongoing research efforts provide meaningful theoretical results \cite{de2000stability} and numerical evidence over numerous successful applications \cite{qin2000overview}.

A quite common approach to nonlinear problems is \emph{explicit MPC} \cite{alessio2009survey}, which relies on a set of precomputed solutions rather than solving onboard the OCP.
However, despite successful applications to aerospace problems are documented in the present literature \cite{hegrenaes2005spacecraft,leomanni2014explicit}, explicit MPC is an intrinsically suboptimal methodology and it may require formidable resources when dealing with high-dimensional OCPs.
Recently,  machine learning techniques have been also proposed to address nonlinear guidance and control problems \cite{izzo2018survey}.
For instance, 
a deep neural network (DNN) can be trained to imitate an expert behavior, as in behavioral cloning \cite{Wang2020asteroid}, or learn an optimal control policy from repeated interactions with the environment, as in reinforcement learning \cite{federicideep2021, zavoli2021reinforcement}. 
However, despite the low evaluation times and the high accuracy in function approximation of DNNs, the absence of theoretical guarantees and the extensive computational effort associated with the training, still limit their application to much easier problems than the one here investigated.



On the other hand, convex optimization techniques are natural candidates for real-time applications, as convex OCPs can be solved by means of highly efficient interior point algorithms, which converge to the optimal solution in polynomial time regardless of the initialization point \cite{boyd2004convex}.
The main challenge consists in formulating a nonlinear (specifically, nonconvex) problem as an equivalent convex one.
Lossless convexification techniques are extremely powerful in this sense, as they can transform problems without introducing any approximations thanks to a convenient change of variables and suitable constraint relaxations that are theoretically guaranteed to retain the same solution as the original problem \cite{accikmecse2011lossless}.
When no lossless convexification exist for the problem at hand, one can always rely on successive convexification methods (first and foremost, successive linearization), which define a sequence of convex problems that eventually converges to the optimal solution of the original problem \cite{liu2017survey}.
Convergence of successive convexification methods is not guaranteed in general, except under appropriate assumptions \cite{liu2014solving,mao2016successive,bonalli2019trajectory}, but extensive numerical evidence suggests the validity of successive convexification in a broad range of applications, 
including 
rocket ascent \cite{benedikter2019convexascent,benedikter2020convex},
spacecraft rendezvous \cite{benedikter2019convexrendezvous},
proximity operations \cite{lu2013autonomous},
low-thrust transfers \cite{wang2018optimization},
rocket landing \cite{acikmese2007convex,szmuk2016successive,sagliano2019generalized},
and atmospheric entry \cite{liu2015entry,wang2017constrained,wang2020improved}.

This paper investigates the use of MPC for the computational guidance of the upper stages of a four-stage VEGA-like vehicle, which features three 
solid rocket motors (SRMs) and a small liquid rocket engine, AVUM, that
performs the final orbit insertion maneuver \cite{vega2014manual}.
%
Due to this stage configuration, this case study poses some unique challenges.
Indeed, the third stage burns out at an extremely high velocity, close to the orbital one, and ends up falling very distant from the launch site.
As a result, the return trajectory of the spent stage must be accurately predicted and constrained to an uninhabited region, which is usually the sea.
To obtain a valuable estimate of the splash-down location, a complete simulation of the return trajectory must be included in the OCP, greatly increasing its complexity.

Uncertainties on the solid rocket motor performance, and especially on its burn time, do not allow to easily pinpoint a splash-down location but rather define a finite-dimension footprint of possible impact points. 
Thus, a further system requirement is bounding the extent of this region.
To this aim, robust optimization and robust MPC can be used to endow the guidance and control with some robustness guarantees.
However, these methods are usually overly conservative since they are primarily based on min-max OCP formulations \cite{bemporad1999robust} or on constraint tightening \cite{chisci2001systems}.
On the other hand, more recent tube-based MPC \cite{langson2004robust} or stochastic optimal control methods, such as chance-constrained optimization \cite{mesbah2016stochastic} or covariance control \cite{chen2016optimal1,chen2016optimal2}, could also be considered as viable options, but their application to the case of non-additive disturbances, such as those that arise from dynamics under uncertain time-lengths, can be quite challenging. 


For all these reasons, VEGA currently relies on a neutral axis maneuver to robustly minimize the footprint size and make the return as insensitive as possible to the third stage performance dispersions \cite{giannini2013vega}.
This maneuver is based on the \emph{null miss condition} developed for ballistic missiles \cite{regan1993dynamics}: over the last few seconds of operation, the stage is constrained to hold on to an attitude such that the splash-down point is retained regardless of any additional velocity increments.
The maneuver reduces the carrying capacity of the launch vehicle since part of the propulsive energy is spent in a non-optimal direction, but it robustly guarantees that the actual splash-down location is sufficiently close to the predicted one, despite additional unpredicted seconds of operation of the SRM.
However, the neutral axis maneuver can be quite difficult to design with traditional methods and usually extensive trajectory validation and verification tasks are required before each launch.
As an alternative to the neutral axis maneuver, the use of retro-rockets to ignite after the third stage separation has been proposed \cite{martens2016innovative}.
However, such a solution would need the integration of additional hardware and mass into a consolidated architecture.

This paper presents an improved version of the MPC algorithm previously proposed by the authors in Ref.~\cite{benedikter2020autonomous}.
The main difference with the past work concerns the introduction of a constraint on the control direction at the stage burnout, which must correspond to the neutral-axis attitude.
In practice, two return trajectories (a nominal and a perturbed one) are considered, whose initial state (i.e., conditions at motor cut-off) differs only by an excess of velocity 
provided in the control direction.
By constraining the impact point of the two return phases to the same value, the neutral axis attitude is indirectly attained.
The envisaged approach overcomes traditional issues in the design of the neutral axis maneuver, which relies on  
cumbersome and time-consuming design procedures, and allows its integration in a modern computational guidance scheme, ensuring the system robustness to extra burn seconds of the SRM.


The paper is organized as follows. 
Section~\ref{sec:problem} describes the guidance and control problem, detailing the phases, the dynamical model, and mission requirements of the OCP to solve.
Section~\ref{sec:cvx} illustrates the state-of-the-art convexification techniques and the $hp$ pseudospectral method used to cast the nonlinear OCP as a sequence of convex problems that rapidly converges to the optimal solution.
In Section~\ref{sec:mpc}, the convex optimization algorithm is embedded into the MPC architecture and a few relaxations are introduced to ensure the recursive feasibility of the algorithm.
Finally, Section~\ref{sec:results} analyzes the effectiveness of the proposed guidance through a Monte Carlo analysis.
In particular, the robustness of the employed algorithm is verified in the presence of 
off-nominal initial conditions,
generic uncertainties on the dynamical model, 
and random SRM performance.
A conclusion section ends the manuscript.

\section{Problem Description}
\label{sec:problem}

MPC is based on solving at every control step an OCP.
In this section, the guidance problem of the third stage of the VEGA-like launch vehicle is described and formulated as a multi-phase OCP.
First, the phase sequence of the OCP is outlined;
then, the considered dynamical model is presented;
and finally, the objective function and all the constraints of the OCP are reported. 

\subsection{Phase Sequence}


\begin{figure}
    \centering
    \includegraphics[width=\columnwidth]{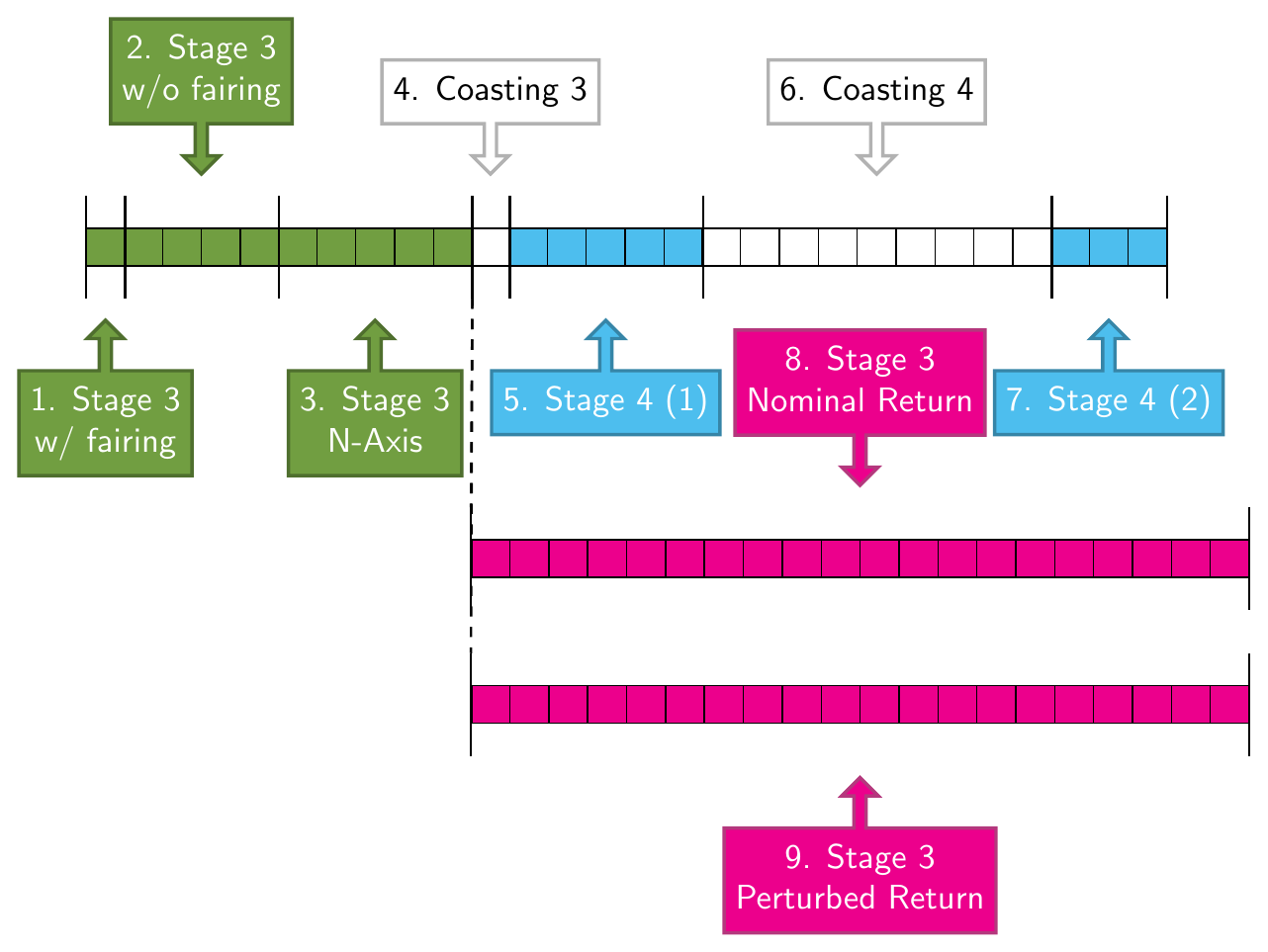}
    \caption{Phases of the optimal control problem.}
    \label{fig:phases}
\end{figure}

Figure~\ref{fig:phases} illustrates schematically the considered phase sequence.
The phases are numbered from \StageThreeOne{} to \StageFourTwo{} in a chronological order, while Phases \Return{} and \PReturn{} simulate the nominal and perturbed uncontrolled return of the spent third stage and start at the burnout of the third stage.
The third stage operation is split into three phases to account for the mass discontinuity associated with the fairing jettisoning, which takes place few seconds after the stage ignition and divides Phase \StageThreeOne{} from Phase \StageThreeTwo{}, and for the transition to the neutral axis maneuver, between Phases \StageThreeTwo{} and \StageThreeThree{}.
In Phase \StageThreeThree{}, the control rate is constrained to a small value, to accomodate the neutral axis maneuver.
The introduction of Phase~\StageThreeThree{} near the end of the stage operation prevents numerical issues to the discretization of rapidly changing dynamics since, during the transition to the neutral axis attitude, the rocket rotates by approximately \SI{40}{\degree} in a few seconds.

A short coasting (Phase \CoastingThree{}) follows the third stage separation, then the fourth stage performs the orbit insertion maneuver, which is composed of two propelled arcs (Phases \StageFourOne{} and \StageFourTwo{}) separated by a long coasting (Phase \CoastingFour{}).
Two return phases are included in the OCP to robustly ensure the splash-down constraint even in the case of additional burn seconds of the SRM.
The initial position and velocity of the nominal return (Phase \Return{}) correspond to the state values at the third stage burnout, while the perturbed return (Phase \PReturn{}) is assumed to start with a larger velocity.
The magnitude of the velocity perturbation is fixed and its value is in the order of the maximum extra $\Delta V$ that can be provided by the SRM, while the direction is constrained to be the same as the control direction at the burnout.
Note that since both returns are constrained to fall at the same splash-down point, the optimal direction of the $\Delta V$ is the neutral axis direction.

Hereinafter, let \smash{$t_0^{(i)}$} and \smash{$t_f^{(i)}$} denote the initial and final time of the $i$-th phase.
For the sake of simplicity, if no phase superscript is specified, then $t_0$ and $t_f$ denote the ignition time of the third stage $\smash{t_0^{(\StageThreeOne{})}}$ and the fourth stage burnout $\smash{t_f^{(\StageFourTwo{})}}$, respectively.
Likewise, let $t_R$ denote the nominal return time of the spent stage \smash{$t_f^{(\Return{})}$} and $t_{PR}$ the perturbed return time \smash{$t_f^{(\PReturn{})}$}.

\subsection{System Dynamics}

The vehicle is modeled as a point mass subject to a 3-DoF translational motion. 
The rotational dynamics are neglected and the rocket axis is assumed to be aligned with the thrust direction at any time. 
The state $\bm{x}$ of the rocket is described by its position $\bm{r}$, velocity $\bm{v}$, and mass $m$, as $\bm{x} = [ x \; y \; z \; v_x \; v_y \; v_z \; m ]$.
The position and velocity are expressed in Cartesian Earth-centered inertial (ECI) coordinates, with the $x$ and $y$ axes in the equatorial plane and the $z$ axis aligned with the Earth's angular velocity $\bm{\omega}_E$, forming a right-hand frame.
The main advantage of using this set of state variables over, for example, a spherical coordinate system \cite{benedikter2021convex}, is the fact that any trajectory can be studied, even missions to polar and high inclination orbits, which are the typical targets for VEGA, without running into singularities.

The only forces assumed to act on the launch vehicle are gravity, the aerodynamic drag, and the engine thrust.
A simple inverse-square gravitational model is assumed, thus the gravity acceleration is $\bm{g} = - \mu \bm{r} / r^3 $,
with $\mu$ denoting the Earth's gravitational parameter.
The drag force is modeled as
$\bm{D} = - 0.5 \rho {v}_{\text{rel}}^2 S C_D \bm{\hat{v}}_{\text{rel}}$,
where 
$\rho$ is the atmospheric density, 
$S$ is the reference surface, 
$C_D$ is the drag coefficient, assumed to be constant, 
and $\bm{v}_{\text{rel}} = \bm{v} - \bm{\omega}_E \times \bm{r}$ is the relative-to-atmosphere velocity.
Lift is neglected as it is much smaller than drag for the system at hand.

The thrust $T$ generated by the SRMs and the associated mass flow rate $\dot{m}_e$ generally vary over time according to the geometry of the solid propellant and, in general, to the design of the stage.
So, every SRM is characterized by a vacuum thrust law $T_{\text{vac}}(t)$ 
and a corresponding mass flow rate law $\dot{m}_e(t)$.
Although upper stages operate at high altitudes, the net thrust is slightly different from the vacuum one due to the the atmospheric pressure $p$ contribution, thus $T = T_{\text{vac}}(t) - p A_e$, where $A_e$ is the nozzle exit area.

Instead, the thrust and mass flow rate of the last stage, which is equipped with a small liquid rocket engine, are assumed to be constant and equal to the maximum attainable values. 
Also, differently from the other stages, AVUM can be cut off and re-ignited. 
So, its burn time can be actively controlled and is determined by the optimization process with the aim of minimizing the consumed propellant mass.

In the propelled phases, the thrust direction $\bm{\hat{T}}$ is a control variable and its components are expressed in the same frame as the velocity.
Being a unit vector, the components satisfy the following identity at any time
\begin{equation}
    \hat{T}_x^2 + \hat{T}_y^2 + \hat{T}_z^2 = 1
	\label{eq:thrust_direction_equality_path_con}
\end{equation}

The resulting equations of motion $\bm{\dot{x}} = \bm{f}(\bm{x}, \bm{u}, t)$ are
\begin{align}
    \dot{\bm{r}} &= \bm{v} \label{eq:original_ODE_r} \\
    \dot{\bm{v}} &= -\frac{\mu}{r^3} \bm{r} + \frac{T}{m} \bm{\hat{T}} - \frac{D}{m} \bm{\hat{v}}_{\text{rel}} \label{eq:original_ODE_v} \\
    \dot{m} &= -\dot{m}_e \label{eq:original_ODE_m}
\end{align}

\subsection{Optimal Control Problem}

The goal of the optimization is to minimize the propellant consumed during the fourth stage operation.
Indeed, if the guidance algorithm minimizes the propellant mass, not only it will provide better responsiveness in off-nominal conditions, but it will also increase the launch vehicle nominal payload, as smaller fuel tanks will lead to an increment of the vehicle carrying capacity.
The objective can be equivalently formulated as maximizing the final mass or minimizing its opposite. Thus, a convenient cost function is
\begin{equation}
    J = -m(\tf)
	\label{eq:objective_m_minimize}
\end{equation}

In an MPC framework, the initial condition of the OCP is updated at every control step with the real-time measurements $\tilde{\bm{x}}_0$, therefore the initial state is completely assigned
\begin{equation}
    \bm{x}(t_0) = \tilde{\bm{x}}_0
    \label{eq:x0}
\end{equation}
The final conditions of Phase \StageFourTwo{} are the following
\begin{align}
    x(\tf)^2 + y(\tf)^2 + z(\tf)^2 &= a_{\text{des}}^2 \label{eq:final_radius_nonlinear} \\
    v_{x}(\tf)^2 + v_{y}(\tf)^2 + v_{z}(\tf)^2 &= v_{\text{des}}^2 \label{eq:final_velocity_nonlinear} \\
    \bm{r}(\tf) \cdot \bm{v}(\tf) &= 0 \label{eq:final_radial_velocity} \\
    x(\tf) v_y(\tf) - y(\tf) v_x(\tf) &= h_{z, \text{des}}
    \label{eq:final_ainc_circular}
\end{align}
In the present study, we consider a circular target orbit with semi-major axis $a_{\text{des}}$, thus Eq.~\eqref{eq:final_velocity_nonlinear} with $v_{\text{des}} = \sqrt{\mu / a_{\text{des}}}$ and Eq.~\eqref{eq:final_radial_velocity} are sufficient to impose the null eccentricity constraint.
As for the inclination, $i_{\text{des}}$ is imposed by prescribing the $z$-component of the angular momentum vector $h_{z, \text{des}} = \cos i_{\text{des}} \sqrt{\mu a_{\text{des}}}$.

To constrain the splash-down of the spent stage, the following conditions are imposed at the end of Phases \Return{} and \PReturn
\begin{align}
    x(t_R)^2 + y(t_R)^2 + z(t_R)^2 &= R_E^2 \label{eq:final_r_reentry_nonlinear} \\
    x(t_{PR})^2 + y(t_{PR})^2 + z(t_{PR})^2 &= R_E^2 \label{eq:final_r_preentry_nonlinear} \\
    z(t_R) = z(t_{PR}) &= z_{R, \text{des}} \label{eq:final_LAT_reentry}    
\end{align}
with $z_{R, \text{des}} = R_E \sin\varphi_{R, \text{des}}$.
Equations~\eqref{eq:final_r_reentry_nonlinear} and \eqref{eq:final_r_preentry_nonlinear} constrain the final radius to be equal to the Earth radius $R_E$.
Instead, Eq.~\eqref{eq:final_LAT_reentry} constrains the latitude of the splash-down point to a desired value $\varphi_{R, \text{des}}$.
For missions toward polar or quasi-polar orbits (e.g., Sun-synchronous orbits), which are the scope of the present paper,
this is equivalent to constraining the distance from the launch site, which is the only quantity that can be actively controlled during the third stage flight.
Indeed, constraining the longitude to an arbitrary value would be very expensive since it would require a maneuver to change the plane of the ascent trajectory.

The phases of the OCP must be linked by some conditions at the interfaces.
In particular, all state variables are continuous across phases, with the relevant exception of the mass, which is discontinuous at the fairing jettisoning and at the third stage separation. 
Thus,
\begin{align}
    m(t_0^{(\StageThreeTwo{})}) &= m(t_f^{(\StageThreeOne{})}) - m_{\text{fairing}} \label{eq:mass_hs_lkg_con} \\
    m(t_0^{(\CoastingThree{})}) &= m(t_f^{(\StageThreeTwo{})}) - m_{\text{dry}, 3} \label{eq:mass_lkg_con_3}
\end{align}
As for Phases \Return{} and \PReturn{}, their initial position and velocity are equal to the burnout conditions of the third stage, with a velocity increment $\bm{\Delta V}_{NA}$ for Phase \PReturn{}, while their mass is equal to the dry mass of the stage:
\begin{align}
    \bm{r}(t_0^{(\Return{})}) &= \bm{r}(t_0^{(\PReturn{})}) = \bm{r}(t_f^{(\StageThreeThree{})}) \label{eq:return_lkg_position} \\ 
    \bm{v}(t_0^{(\Return{})}) &= \bm{v}(t_f^{(\StageThreeThree{})}) \label{eq:return_lkg_velocity} \\ 
    \bm{v}(t_0^{(\PReturn{})}) &= \bm{v}(t_f^{(\StageThreeThree{})}) + \bm{\Delta V}_{NA} \label{eq:preturn_lkg_velocity} \\ 
    m(t_0^{(\Return{})}) &= m(t_0^{(\PReturn{})}) = m_{\text{dry}, 3} \label{eq:return_mass}
\end{align}

We also impose a continuity constraint between the velocity perturbation $\bm{\Delta V}_{NA}$ and the final control direction at the nominal burnout of the third stage. Thus,
\begin{equation}
    \bm{\Delta V}_{NA} / \Delta V_{NA} = \bm{\hat{T}}(t_f^{(\StageThreeThree{})})
    \label{eq:pert_reentry_DV_in_dir3f}
\end{equation}
This ensures that in case of an SRM longer burn time, the launch vehicle would hold on to this direction, thus minimizing the shift of the impact point.

Finally, since the fairing is jettisoned as soon as possible to reduce the dry mass of the system, from that moment on, the thermal flux on the rocket must be lower than a given value to preserve the integrity of the payload.
Thus, a constraint on the maximum heat flux is included in the formulation for Phases \StageThreeTwo{} to \StageFourTwo{}
\begin{equation}
    \dot{Q} = \frac{1}{2} \rho v_{\text{rel}}^3 \leq \dot{Q}_{\text{max}}
    \label{eq:heat_flux_nonlinear}
\end{equation}

\section{Convex Formulation}
\label{sec:cvx}

This section outlines the convexification strategy devised to convert the OCP into a sequence of convex problems that quickly converges to the optimal solution.
This is a crucial step in the design of the guidance algorithm, as it greatly reduces the cost of the onboard optimization and enables the MPC controller to attain high update frequencies.



First, lossless convexification strategies, such as changes of variables and constraint relaxations, are employed to reduce the nonconvexity of the original problem.
Instead, the remaining nonconvexities are tackled via successive linearization.
Safeguarding modifications, such as virtual controls and a trust region on the time-lengths of the phases, are also included in the formulation.
The resulting convex OCP is discretized via a $hp$ Radau pseudospectral method to set up a finite-dimensional optimization problem.
Finally, the reference solution update strategy and the convergence criteria are outlined. 

\subsection{Change of Variables}

In Eqs.~\eqref{eq:original_ODE_r}--\eqref{eq:original_ODE_m}, the control variable $\bm{\hat{T}}$ is coupled with the state variable $m$.
This may lead to numerical issues, such as high-frequency jitters in the control law \cite{liu2015entry}.
Instead, the problem should be formulated in a way such that the dynamics are linear in the controls, for instance
\begin{equation}
    \bm{f} = \tilde{\bm{f}}(\bm{x}, t) + \tilde{B} \bm{u}
    \label{eq:ODE_split}
\end{equation}

This can be achieved by introducing a new control
\begin{equation}
    \bm{u} = \frac{T}{m} \bm{\hat{T}} \label{eq:control_vector_cvx}
\end{equation}
which includes both the thrust-to-mass ratio $T / m$ and the thrust direction vector $\bm{\hat{T}}$.
By replacing it in Eqs.~\eqref{eq:original_ODE_r}--\eqref{eq:original_ODE_m}, control-affine dynamics are obtained
\begin{align}
    \dot{\bm{r}} &= \bm{v} \label{eq:affine_ODE_r} \\
    \dot{\bm{v}} &= -\frac{\mu}{r^3} \bm{r} + \bm{u} 
    - \frac{D}{m} \bm{\hat{v}}_{\text{rel}} 
    \label{eq:affine_ODE_v} \\
    \dot{m} &= -\dot{m}_e \label{eq:affine_ODE_m}
\end{align}
and the $\tilde{B}$ matrix in Eq.~\eqref{eq:ODE_split} is
$\tilde{B} = [\bm{0}_{3 \times 3} \;
    \bm{I}_{3 \times 3} \;
    \bm{0}_{1 \times 3}]^T$
with $\bm{0}_{m \times n}$ and $\bm{I}_{m \times n}$ denoting the null and identity matrix of size $m \times n$.


Also the new control variables must satisfy Eq.~\eqref{eq:thrust_direction_equality_path_con}, which is reformulated as
\begin{equation}
    u_x^2 + u_y^2 + u_z^2 = u_N^2
	\label{eq:thrust_direction_equality_path_con_new}
\end{equation}
Note that we introduced an additional variable $u_N$, constrained as
\begin{equation}
    u_N = \frac{T}{m} \label{eq:u_N}
\end{equation}
The introduction of this nonlinear path constraint is the price to pay to obtain control-affine dynamics.
Nevertheless, the advantages of this change of variables outweigh the disadvantages, since, despite the inclusion of the additional constraint, the complexity of the problem is reduced compared to a coupling of state and control variables in the dynamics.

\subsection{Constraint Relaxation}
\label{subsec:constraint_relaxation}
Equation~\eqref{eq:thrust_direction_equality_path_con_new} is a nonlinear equality constraint that must be convexified.
This constraint belongs to a class of constraints that is suitable for a lossless relaxation \cite{acikmese2007convex,accikmecse2011lossless,liu2016exact,benedikter2020convex}.
So, by replacing the equality with an inequality, a second-order cone constraint is formulated
\begin{equation}
    u_x^2 + u_y^2 + u_z^2 \leq u_N^2
	\label{eq:thrust_direction_cone_con}
\end{equation}
This relaxation is lossless because, despite it defines a larger feasible set, the solution of the convex problem is the same as the original and Eq.~\eqref{eq:thrust_direction_cone_con} is satisfied with the equality sign.
A theoretical proof of the lossless property of this relaxation can be found in Proposition~1 of Ref.~\cite{benedikter2020convex}.
It should be noted that Proposition~1 holds as long as the heat flux path constraint \eqref{eq:heat_flux_nonlinear} is inactive almost 
everywhere (i.e., it is active, at most, at isolated points in time).
However, this is a mild assumption for the problem at hand, as the heat flux profile is rarely active over arcs of finite duration.


\subsection{Successive Linearization}
The remaining nonlinear expressions are replaced with first-order Taylor series expansions around a reference solution that is updated at every iteration.

\subsubsection{Equations of Motion}
The equations of motion \eqref{eq:affine_ODE_r}--\eqref{eq:affine_ODE_m} are affine in the control but still nonlinear in the state; thus, they must be linearized.
First, we exploit a common strategy in direct optimization methods to account for free-time phases \cite{betts2010practicalDt}.
We introduce an independent variable transformation from physical time $t$ to a new variable $\tau$ defined, for each phase, over a fixed unitary domain $[0, 1]$.
The relationship between the two independent variables is the following
\begin{equation}
    t^{(i)} = t_0^{(i)} + (t_f^{(i)} - t_0^{(i)}) \, \tau
\end{equation}
Note that the time dilation $s$ between $t$ and $\tau$ corresponds to the time-length of the phase
\begin{equation}
    s^{(i)} = \frac{d t}{d \tau}  = t_f^{(i)} - t_0^{(i)}
    \label{eq:time_dilation}    
\end{equation}
The time-length $s^{(i)}$ is then introduced as an additional optimization variable for every free-time phase.

The dynamics are thus expressed in terms of $\tau$ and linearized around a reference solution $\{\bar{\bm{x}}, \bar{\bm{u}}, \bar{s}\}$
\begin{equation}
    \bm{x}' \vcentcolon= \frac{d \bm{x}}{d \tau} = s \bm{f}\left( \bm{x}, \bm{u}, \tau \right) \approx A \bm{x} + B \bm{u} + \Sigma s + \bm{c}
    \label{eq:linear_ODEs}
\end{equation}
where the following matrices and vectors were introduced
\begin{align}
    A &= \bar{s} \frac{\partial \bm{f}}{\partial \bm{x}} ( \bar{\bm{x}}, \bar{\bm{u}}, \tau ) \label{eq:A_matrix_definition} \\
    B &= \bar{s} \frac{\partial \bm{f}}{\partial \bm{u}} ( \bar{\bm{x}}, \bar{\bm{u}}, \tau ) \label{eq:B_matrix_definition} \\
    \Sigma &= \bm{f} ( \bar{\bm{x}}, \bar{\bm{u}}, \tau ) \label{eq:P_matrix_definition} \\
    \bm{c} &= - (A \bar{\bm{x}} + B \bar{\bm{u}}) \label{eq:C_vector_definition}
\end{align}

Finally, since the linearization may lead to artificial infeasibility \cite{mao2016successive}, a virtual control signal $\bm{q}$ is added to the dynamics
\begin{equation}
    \bm{x}' = A \bm{x} + B \bm{u} + \Sigma s + \bm{c} + \bm{q}
    \label{eq:linear_ODEs_vc}
\end{equation}
This guarantees the feasibility of the convex problem.
Naturally, as the virtual controls are unphysical variables, their use must be highly penalized.
Thus, the following penalty term is added to the cost function
\begin{equation}
    J_q = \lambda_q P(\bm{q})
    \label{eq:penalty_vc}
\end{equation}
where $\lambda_q$ is a sufficiently large penalty weight and $P(\bm{q})$ a suitable penalty function, which will be defined after discretization in Section~\ref{subsec:discretization}.

\subsubsection{Boundary Constraints}
The final conditions at the payload release \eqref{eq:final_radius_nonlinear}--\eqref{eq:final_r_reentry_nonlinear} are linearized as
\begin{align}
    \bar{\bm{r}}(t_f) \cdot \bar{\bm{r}}(t_f) 
    + 2 \bar{\bm{r}}(t_f) \cdot (\bm{r}(t_f) - \bar{\bm{r}}(t_f)) 
    &= a_{\text{des}}^2 
    \label{eq:final_radius_linearized} \\
    \bar{\bm{v}}(t_f) \cdot \bar{\bm{v}}(t_f) 
    + 2 \bar{\bm{v}}(t_f) \cdot (\bm{v}(t_f) - \bar{\bm{v}}(t_f)) 
    &= v_{\text{des}}^2 
    \label{eq:final_velocity_linearized} \\
    \bar{\bm{r}}(t_f) \cdot \bar{\bm{v}}(t_f) 
    + \bar{\bm{v}}(t_f) \cdot (\bm{r}(t_f) - \bar{\bm{r}}(t_f))&
    \nonumber \\
    + \bar{\bm{r}}(t_f) \cdot (\bm{v}(t_f) - \bar{\bm{v}}(t_f)) &= 0
    \label{eq:final_radial_velocity_linearized} \\
    \bar{v}_y(t_f) (x(t_f) - \bar{x}(t_f))
    - \bar{v}_x(t_f) (y(t_f) - &\bar{y}(t_f))
    \nonumber \\
    - \bar{y}(t_f) v_x(t_f) + \bar{x}(t_f) v_y(t_f) &= h_{z, \text{des}}
    \label{eq:final_ainc_linearized_circular}
\end{align}
Likewise, the terminal radius constraints at the splash-down, Eqs.~\eqref{eq:final_r_reentry_nonlinear} and \eqref{eq:final_r_preentry_nonlinear}, are formulated as
\begin{equation}
    \bar{\bm{r}}(t_f^{(i)}) \cdot \bar{\bm{r}}(t_f^{(i)}) 
    + 2 \bar{\bm{r}}(t_f^{(i)}) \cdot (\bm{r}(t_f^{(i)}) - \bar{\bm{r}}(t_f^{(i)})) 
    = R_E^2
    \label{eq:final_r_reentry_linearized}
\end{equation}
with $i = \Return{}, \PReturn{}$.

The linearization of the above constraints may cause artificial infeasibility. 
Thus, virtual buffer zones are introduced.
In particular, Eqs.~\eqref{eq:final_radius_linearized}--\eqref{eq:final_r_reentry_linearized} are grouped into a vector $\bm{\chi} = \bm{0}$ that is then relaxed as $\bm{\chi} = \bm{w}$.
The vector $\bm{w}$ holds all the virtual buffers, which are highly penalized by adding the following penalty term to the cost function
\begin{equation}
    J_w = \lambda_w \norm{\bm{w}}_1
\end{equation}
with $\lambda_w$ denoting the penalty weight of the virtual buffers. 


\subsubsection{Path Constraints}
The heat flux constraint \eqref{eq:heat_flux_nonlinear} is another nonlinear expression that must be linearized.
By considering its Taylor series expansion, we obtain
\begin{equation}
    \frac{\dot{\bar{Q}} + \frac{\partial \dot{\bar{Q}}}{\partial \bm{r}} \cdot (\bm{r} - \bar{\bm{r}}) + \frac{\partial \dot{\bar{Q}}}{\partial \bm{v}} \cdot (\bm{v} - \bar{\bm{v}})}{\dot{Q}_{\text{max}} } \leq 1
    \label{eq:heat_flux_linearized}
\end{equation}
where
\begin{align}
    \frac{\partial\dot{\bar{Q}}}{\partial \bm{r}} &= \frac{1}{2} \frac{d \bar{\rho}}{d \bm{r}} \bar{v}_{\text{rel}}^3 + \frac{3}{2} \bar{\rho} \bar{v}_{\text{rel}} \bm{\omega}_E \times \bm{\bar{v}}_{\text{rel}} \label{eq:dQ_dr} \\
    \frac{\partial\dot{\bar{Q}}}{\partial \bm{v}} &= \frac{3}{2} \bar{\rho} \bar{v}_{\text{rel}} \bm{\bar{v}}_{\text{rel}} \label{eq:dQ_dv}
\end{align}
Note that to avoid numerical issues, we normalized the constraint with respect to the maximum heat flux.

Likewise, the auxiliary variable $u_N$ must satisfy Eq.~\eqref{eq:u_N} at every time.
This is a nonlinear path constraint and it is linearized as
\begin{equation}
    u_N = \frac{T_{vac} - \bar{p} A_e}{\bar{m}} \left(2 - \frac{m}{\bar{m}} \right) - 
    \frac{A_e}{\bar{m}} \frac{d \bar{p}}{d \bm{r}} \cdot (\bm{r} - \bar{\bm{r}})     
    \label{eq:u_N_path_con_linearized}
\end{equation}

\subsection{Trust Region on Time-Lengths}
Thanks to the change of variables \eqref{eq:control_vector_cvx},
$\bm{f}$ is linear in the control and the $A$ and $B$ matrices do not depend on $\bar{\bm{u}}$.
However, $\Sigma$ and $\bm{c}$ are functions of the reference controls and the successive convexification algorithm may suffer from instability issues, as intermediate controls can change significantly in the first iterations \cite{benedikter2019convexrendezvous}.
Nevertheless, when $s \approx \bar{s}$, Eq.~\eqref{eq:linear_ODEs} reduces to
\begin{equation}
    \bm{x}' = A \bm{x} + B \bm{u} + \tilde{\bm{c}}
    \label{eq:linear_ODEs_t_fixed}
\end{equation}
with $\tilde{\bm{c}} = \bar{s} \tilde{\bm{f}}(\bar{\bm{x}}, \tau) - A \bar{\bm{x}}$, 
thus entirely dropping the dependence on $\bar{\bm{u}}$.

To prevent excessive changes from the reference value, a soft trust region constraint is imposed on the time-lengths $s$ 
\begin{equation}
    | s^{(i)} - \bar{s}^{(i)} | \leq \delta^{(i)} \qquad i = \CoastingFour{}, \StageFourTwo{}
    \label{eq:Dt_trust_region}
\end{equation}
This constraint is imposed only on the duration of Phases \CoastingFour{} and \StageFourTwo{} since only these arcs resulted particularly sensitive to instability phenomena \cite{benedikter2020convex}. 
The trust radii $\delta^{(i)}$ are additional optimization variables bounded in a fixed interval $[0, \smash{\delta_{\text{max}}^{(i)}}]$ and slightly penalized in the cost function via the following penalty terms
\begin{equation}
    J_\delta^{(i)} = \lambda_{\delta}^{(i)} \delta^{(i)}
    \qquad i = \CoastingFour{}, \StageFourTwo{}
    \label{eq:delta_penalty_term}
\end{equation}
In the authors' experience, a suitable choice of the upper bound \smash{$\delta_{\text{max}}^{(i)}$} is usually somewhere between 1\% and 10\% of $\bar{s}^{(i)}$.
Instead, the penalty weights $\lambda_{\delta}$ should be picked as small as possible in order not to converge toward suboptimal solutions (e.g., in the range $\smash{10^{-6}} \div \smash{10^{-3}}$).

\subsection{Discretization}
\label{subsec:discretization}

The problem formulated so far is infinite-dimensional since state and controls are continuous-time functions.
However, to solve it numerically, the problem must be cast as a finite set of variables and constraints.
In this paper, we use a $hp$ pseudospectral method, which is a particularly performing transcription method since it combines the advantages of $h$ and $p$ schemes.
Indeed, pseudospectral methods exhibit exponential convergence rate only in regions where the solution is smooth, but combining them with $h$ methods allows to introduce mesh nodes near potential discontinuities, thus extending their appealing properties to problems with rapidly changing dynamics \cite{darby2011hp}.
These features are key to reduce the computational cost of the onboard optimization, as such a discretization method limits the problem dimension and yet it accurately represents the continuous-time problem.

In a $hp$ scheme, the time domain of each phase is split into $h$ subintervals, or \emph{segments}, and the differential constraints are enforced in each segment via local orthogonal collocation. 
The order $p$ of the collocation can vary among the segments, but, for the sake of simplicity, we adopt the same order over all the segments of a phase.
The Radau pseudospectral method (RPM) \cite{garg2011advances} is used as collocation scheme since it is one of the most accurate and performing pseudospectral methods \cite{garg2010unified}.
Furthermore, the Legendre-Gauss-Radau (LGR) points avoid redundant control variables at the segment boundaries and approximate the optimal control at each mesh point, except for the final node of the final segment that must be extrapolated \emph{a posteriori} from the polynomial approximation of the control signal.

Once a grid is defined, the state and control are discretized over it, generating a finite set of variables $(\bm{x}_j^k, \bm{u}_j^k)$.
The superscript $k$ denotes the $k$-th segment, while the subscript $j$ refers to the $j$-th node of the segment.
In each segment $[\tau_{k}, \tau_{k+1}]$, the state and control signals are approximated via Lagrange polynomial interpolation and the derivative of the state approximation is imposed to be equal to the right-hand side of the equations of motion \eqref{eq:linear_ODEs_vc} at the LGR collocation points
\begin{equation}
    \sum_{j = 1}^{p + 1} D_{ij}^k \bm{x}_j^k = \frac{\tau_{k+1} - \tau_{k}}{2} \bm{f}_i^k 
    \qquad 
    i = 1, \dots, p 
    \label{eq:collocation_radau}
\end{equation}
where $\bm{f}_i^k = A_i^k \bm{x}_i^k + B_i^k \bm{u}_i^k + \Sigma_i^k s + \bm{c}_i^k + \bm{q}_i^k$.
The LGR differentiation matrix $D^k$ can be efficiently computed via barycentric Lagrange interpolation \cite{berrut2004barycentric}.

Note that also the virtual control $\bm{q}$ and the linearization matrices \eqref{eq:A_matrix_definition}--\eqref{eq:C_vector_definition} are discretized over the mesh.
We now define also the penalty function $P(\bm{q})$ from Eq.~\eqref{eq:penalty_vc} as the 1-norm of the vector containing all the discrete virtual control variables $\bm{q}_i^k$.
Likewise, path constraints \eqref{eq:heat_flux_linearized} and \eqref{eq:u_N_path_con_linearized} are imposed at every node, generating a finite set of algebraic constraints.
Further details on the implementation of the $hp$ Radau pseudospectral method can be found in Ref.~\cite{patterson2014gpopsii}.


\subsection{Update of the Reference Solution}

Successive linearization is based on Taylor series expansions around a reference solution.
After each convex problem is solved, the reference solution must be updated somehow.
The most common update method simply consists of replacing the reference solution with the new one.
However, to provide further algorithmic robustness to the successive convexification approach, we employ an alternative approach, referred to as \emph{filtering}, that was successfully used by the authors in previous works \cite{benedikter2019convexrendezvous, benedikter2020convex, benedikter2020autonomous}.

Filtering consists of updating the reference solution with a weighted sum of the previous solutions.
By averaging out the values over different solutions, the sequence of reference solutions is smoother and instabilities due to possible diverging intermediate iterations are filtered out.
As a result, the algorithm is less sensitive to artificial unboundedness and trust regions \cite{mao2016successive} or other expedients \cite{wang2020improved}, typically necessary in such a complex problem, can be avoided. 
Indeed, the major appeal of filtering is the fact that it does not introduce any additional constraints or penalty terms to the problem, retaining its original formulation.

In practice, the reference solution for the $i$-th problem is computed as
\begin{equation}
	\bar{x}^{(i)} = \sum_{k = 1}^{K} \alpha_k x^{\max \{0, (i - k)\} }
\end{equation}
where $K$ is the number of previous solutions to account for, $\alpha_k$ are the corresponding weights, 
and $x^{(i)}$ denotes the solution to the $i$-th problem.
Note that until $i \geq K$ the initial reference solution $x^{(0)}$ appears multiple times in the sum.
Three previous solutions are used for the reference solution update, weighted in the same way as in the authors' previous works \cite{benedikter2019convexrendezvous,benedikter2020convex,benedikter2020autonomous}.
The values of the weights are
$\alpha_1 = 6/11$, $\alpha_2 = 3/11$, and $\alpha_3 = 2/11$.
        
Eventually, the sequential algorithm terminates when all the following criteria are met:
\begin{enumerate}[label=(\roman*)]
\item the difference between the computed solution and the reference one converges below an assigned tolerance $\epsilon_{\text{tol}}$
\begin{equation}
    \norm{\bm{x} - \bar{\bm{x}}}_\infty < \epsilon_{\text{tol}}
    \label{eq:successive_cvx_termination_condition}
\end{equation}
\item the computed solution adheres to the nonlinear dynamics within a tolerance $\epsilon_f$
\begin{equation}
    \Big\lVert \sum_{j = 1}^{p + 1} D_{ij}^k \bm{x}_j^k - \frac{\tau_{k+1} - \tau_{k}}{2} \bm{f}_i^k \Big\rVert_\infty < \epsilon_f 
    \label{eq:acceptable_dyn}
\end{equation}
in each phase, for $i = 1, \dots, p$, and $k = 1, \dots, h$; 
\item the virtual buffers of the computed solution are below the dynamics tolerance $\epsilon_{f}$
\begin{equation}
    \norm{\bm{w}}_\infty < \epsilon_{f}
    \label{eq:acceptable_vb}
\end{equation}
\end{enumerate}

 
\section{Model Predictive Control}
\label{sec:mpc}

\subsection{Guidance Strategy}

In the present work, the MPC controller is applied only to the third stage operation, thus the control cycle repeats until the stage burnout, which is the MPC stopping criterion.
However, the \emph{prediction} horizon (i.e., the time domain of the OCP) must last until the first point in time that allows to evaluate performance and account for all the mission constraints, that is, the payload release into orbit and the splash-down of the third stage.


The peculiarity of the present application is a possible mismatch between the expected and the actual burnout time of the stage.
To design a guidance algorithm that is robust to this uncertainty, we devised a strategy that consists of a closed-loop guidance until the minimum burn time and an open-loop neutral axis maneuver in the extra burn seconds.
The neutral axis direction is determined onboard by the OCP, which provides the optimal direction of an additional $\Delta V$ such that the return point is unchanged.



\subsection{Noise Model}
\label{subsec:noise}

To simulate the actual stage operation, a numerical integration of the original dynamical model, Eqs.~\eqref{eq:original_ODE_r}--\eqref{eq:original_ODE_m}, is carried out.
To account for model uncertainty and external disturbance, white Gaussian noise is added to the dynamics, resulting in the following stochastic differential equation (SDE)
\begin{equation}
    d\bm{x} = \bm{f}(\bm{x}, \bm{u}, t) dt + G(\bm{x}, \bm{u}, t) d\bm{w}
    \label{eq:sde}
\end{equation}
where $\bm{w}$ is a $n_w$-dimensional Wiener process and the matrix $G$ determines how the external disturbances affect the system.

In the present case, we assume that the standard deviation of the random disturbance is proportional to the thrust level by a factor $\alpha_G$.
Also, we model two different noise contributions: the first one is aligned with the thrust direction while the other one is in a random direction.
We suppose that the former is $\alpha_{G_\parallel}$ times larger than the latter, as most of the noise is usually in the thrust direction. 
The resulting $G$ matrix is
\begin{equation}
    G = \alpha_G \frac{T}{m} \begin{bmatrix}
        \bm{0}_{3 \times 1} & \bm{0}_{3 \times 3} \\
        \alpha_{G_\parallel} \bm{\hat{T}} & \bm{I}_{3 \times 3} \\
        \bm{0}_{1 \times 1} & \bm{0}_{1 \times 3} 
    \end{bmatrix}
    \label{eq:G_matrix}
\end{equation}
Note that $n_w = 4$.

Also the initial conditions at the stage ignition are perturbed to account for possible performance dispersions of the previous stages.
The position error is modeled as a perturbation of the altitude only, as perturbations in other directions do not significantly affect performance.
Errors on $r$ are uniformly sampled in the range $[-\Delta_r, \Delta_r]$.
Instead, the scattered initial conditions on the velocity are generated by uniformly sampling noise from a sphere of radius $\Delta_v$.

\begin{figure}[h]
    \centering
    \includegraphics[width=\columnwidth]{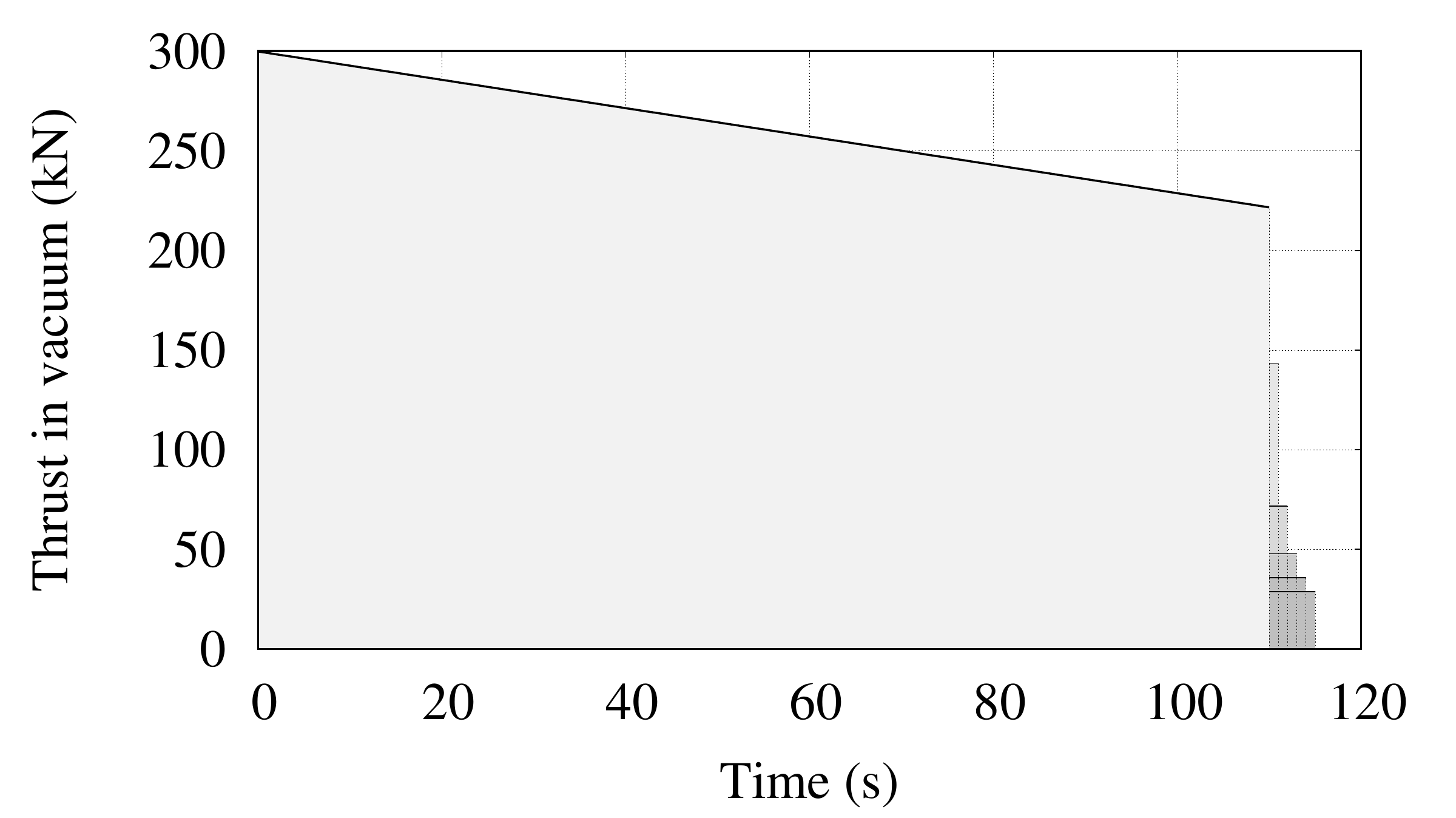}
    \caption{Nominal and perturbed thrust profiles}
    \label{fig:thrust_profiles}
\end{figure}

As for the SRM thrust profile used in the simulation, every simulation considers a perturbed thrust law that features an additional burn time $\Delta t_b$ of up to \SI{5}{\s} and a small increase in the engine total impulse $\Delta I_{\text{tot}}$ up to \num{0.5}\% with respect to the nominal one used as a reference in the OCP.
By randomly scattering these quantities, a constant, nonzero thrust level is assumed in the extra seconds of operation
\begin{equation}
    T_{\text{vac}} (t > t_{b, \text{nom}}) = \frac{\Delta I_{\text{tot}}}{\Delta t_b}
\end{equation}
The nominal thrust profiles and some examples of perturbed laws are shown in Fig.~\ref{fig:thrust_profiles}.

\subsection{Relaxation}

Recursive feasibility of the OCP is key to the design of an MPC controller, as the optimization procedure must provide a control signal at every control step.
The convex problem formulated in Section~\ref{sec:cvx} already includes virtual controls and buffer zones to ensure the feasibility of each instance.
Under nominal circumstances, these expedients would guarantee the feasibility of the OCP.
However, when the system operates in off-nominal conditions, satisfying some constraints may be impossible, and, as a result, the optimization may provide unphysical solutions that actively exploit virtual variables or no solution at all. 
A common strategy to avoid this undesired phenomenon is relaxing the initial condition constraint \eqref{eq:x0}, e.g., as proposed in Ref.~\cite{wang2019optimal}.
This is a general approach that can be readily employed regardless of the application.
However, as such, it does not exploit any knowledge on the problem specificities.
Rather, if possible, it is more efficient to identify the critical mission requirements and formulate a properly relaxed problem.

For the problem under investigation, only the constraint on maximum thermal flux \eqref{eq:heat_flux_linearized} may cause an infeasible problem instance.
Indeed, the nominal trajectory is such that the heat flux is greater than or equal to the threshold just until the fairing jettisoning.
So, due to off-nominal initial conditions or external disturbances in the dynamics, the heat flux may be greater than expected at some points in the simulation and keeping it below threshold may be impossible.
The infeasibility is solved by relaxing the heat flux constraint \eqref{eq:heat_flux_linearized} as
\begin{equation}
    \frac{\dot{\bar{Q}} + \frac{\partial \dot{\bar{Q}}}{\partial \bm{r}} \cdot (\bm{r} - \bar{\bm{r}}) + \frac{\partial \dot{\bar{Q}}}{\partial \bm{v}} \cdot (\bm{v} - \bar{\bm{v}})}{\dot{Q}_{\text{max}} } \leq 1 + \delta_{\dot{Q}}
    \label{eq:heat_flux_linearized_relax}
\end{equation}
where $\delta_{\dot{Q}}$ is a non-negative optimization variable that is penalized in the cost function via the term 
\begin{equation}
    J_{\dot{Q}} = \lambda_{\dot{Q}} \delta_{\dot{Q}}
\end{equation}
In the authors' experience, high values should not be assigned to the penalty weight \smash{$\lambda_{\dot{Q}}$}, since the optimization would then exploit virtual controls and may never meet the converge criterion \eqref{eq:acceptable_dyn}.
Rather, 
considering a conservative nominal threshold \smash{$\dot{Q}_{\text{max}}$} and assigning a small penalty to its violation appears as a much more effective strategy.



\subsection{Update}

At the end of each control cycle, the OCP must be updated with the measurements coming from the navigation system.
The update concerns the initial condition of the OCP, Eq.~\eqref{eq:x0}, which is updated with the state vector resulting from the numerical integration of the SDE~\eqref{eq:sde}.
Also, the reference solution $\{\bar{\bm{x}}, \bar{\bm{u}}, \bar{s}\}$ is replaced with the optimal solution computed at the previous time, removing the portion of flight elapsed between the two control steps.

Since the time-length of the prediction horizon reduces at each update due to the receding-horizon implementation of the MPC algorithm, to save computational resources, the size of the discretization grid is also reduced linearly at every step from the values in Table~\ref{tab:mesh} until a minimum value of \num{5} nodes per segment.
Also, Phases \StageThreeOne{} and \StageThreeTwo{} are removed from the OCP as soon as they are over.

\begin{table}[h]
    \centering
    \caption{Discretization segments and order in each phase}
    \label{tab:mesh}
    \medskip
    \begin{tabular}{l c c c c c c c c c}
    \hline
    \bf Phase & \StageThreeOne{} & \StageThreeTwo{} & \StageThreeThree{} & \CoastingThree{} & \StageFourOne{} & \CoastingFour{} & \StageFourTwo{} & \Return{} & \PReturn{} \\
    \hline
    $h$ & 1 & 1 & 1 & 1 & 1 & 1 & 1 & 5 & 5 \\
    $p$ & 5 & 19 & 5 & 9 & 19 & 19 & 19 & 20 & 20 \\
    \hline
    \end{tabular}
\end{table}

As soon as the nominal burn time $t_{b, \text{nom}}$ is passed, the guidance switches from a closed-loop MPC architecture to an open-loop neutral axis maneuver for the extra seconds of stage operation.
The control direction during the neutral axis maneuver is constant and corresponds to the terminal control action computed in the last MPC cycle.

It is worthwhile mentioning that, in the present work, the time to solve the optimization problem is supposed null, but actually the simulation should introduce a delay between the measurement of the updated state and the actuation of the optimal control law.
This delay can be quite long and introduce significant deviations from the predicted trajectory if the computation time is long.
Nevertheless, thanks to the convexification strategy, the time required to solve the OCP is assumed to be sufficiently brief. 
Preliminary numerical results confirm the validity of this hypothesis.

\section{Numerical Results}
\label{sec:results}

In this section, numerical results are presented to assess the performance and robustness of the devised MPC framework.
In particular, we compare the devised guidance strategy with the one previously proposed by the authors \cite{benedikter2020autonomous}.
We refer to the latter as \emph{single-return} algorithm and to the novel one as \emph{multi-return} since it includes an additional perturbed return phase in the OCP to autonomously determine the neutral axis direction.
The algorithms were implemented in C++ and use Gurobi's \cite{gurobi} second-order cone programming solver to solve the OCP.
All the simulations were run on a computer equipped with Intel\textregistered{} Core\texttrademark{} i7-9700K CPU @ \SI{3.60}{\giga\Hz}.

The parameters of the successive convexification algorithm are the same as in Ref.~\cite{benedikter2020convex}.
Specifically, the penalty weights on the time-lengths, $\smash{\lambda_\delta^{(\CoastingFour{})}}$ and $\smash{\lambda_\delta^{(\StageFourTwo{})}}$, were set to $\smash{\num{e-4}}$, while, to highly penalize solutions that actively exploit virtual variables, we set \smash{$\lambda_q = \lambda_w = \num{e4}$}.
The tolerances on the convergence criteria, Eqs.~\eqref{eq:successive_cvx_termination_condition}--\eqref{eq:acceptable_vb}, were prescribed as \smash{$\epsilon_{\text{tol}} = \num{e-4}$} and \smash{$\epsilon_f = \num{e-6}$}.

\subsection{Nominal Trajectories}

The investigated scenario is the same as in Refs.~\cite{benedikter2020convex} and \cite{benedikter2020autonomous},
that is, an ascent trajectory from an equatorial launch base to a \SI{700}{\kilo\m} circular orbit with $i_{\text{des}} = \SI{90}{\degree}$ and a splash-down of the third stage constrained to $\varphi_{R, \text{des}} = \SI{60}{\degree}$.
Also, the same data are used to model the launch vehicle stages and aerodynamic coefficients, with the exception of the uncertain additional burn time of the third stage, which is modeled as detailed in Section~\ref{sec:mpc}.
The U.S. Standard Atmosphere 1976 model is used to evaluate the air density and pressure as functions of the altitude~\cite{us1976atm}.

The third stage nominal burn time is $t_{b, 3} = \SI{110}{\s}$, which is also the minimum burn time considered in the simulations.
The fairing is released $\Delta t^{(\StageThreeOne{})} = \SI{5.4}{\s}$ after ignition, while Phase \StageThreeThree{} is prescribed to last $\Delta t^{(\StageThreeThree{})} = \SI{5}{\s}$.
Thus, Phase \StageThreeTwo{} lasts $\Delta t^{(\StageThreeTwo{})} = \SI{99.6}{\s}$.
The duration of the coasting arc after the third stage separation (Phase \CoastingThree{}) is fixed to $\Delta t^{(\CoastingThree{})} = \SI{15.4}{\s}$.
The time-lengths of all the following phases are left to be optimized.

In the OCP, the threshold on the heat flux used to compute the solution is set to $\dot{Q}_{\text{max}}=\SI{900}{\W\per\square\m}$.
Since $\dot{Q}_{\text{max}}$ is lower than the maximum bearable heat flux that must not be exceeded in simulations (\SI{1135}{\watt\per\square\m}) 
\cite{vega2014manual}, the relaxation of the heat flux constraint with a relatively small value of the penalty coefficient (\smash{$\lambda_{\dot{Q}} = \num{e-2}$}) is justified.


The nominal trajectories of the single and multi-return algorithms were computed via the convex optimization approach described in Ref.~\cite{benedikter2020convex}, with the only difference of including Phases \StageThreeThree{} and \PReturn{} in the OCP for the multi-return scenario.
The considered velocity increment at the start of Phase \PReturn{} is $\Delta V_{NA} = \SI{25}{\m\per\s}$, which is comparable to the (expected) maximum velocity increment provided by the SRM compared to the nominal performance.
The nominal trajectories are used as initial reference solution in the first control cycle and provide the nominal payload masses: \SI{1400.1}{\kilo\g} for the single-return problem and \SI{1328.3}{\kilo\g} for the multi-return one. 
The carrying capacity of the multi-return strategy is smaller due to the fact that the third stage must be oriented in the neutral axis direction at burnout, while the single-return problem does not take into account such requirement.
Thus, the difference in the payload mass represents the price to pay for attaining the neutral axis attitude.

\subsection{Monte Carlo Campaigns}

A Monte Carlo analysis on the combined effect of off-nominal initial conditions, in-flight disturbance, and uncertain thrust profiles as described in Section~\ref{subsec:noise} was carried out for both the single-return strategy and the present one.
The initial radius error is scattered uniformly in a range of size $\Delta_r = \SI{500}{\m}$ and the initial velocity error is sampled from a sphere of radius $\Delta_v = \SI{40}{\m\per\s}$.
Three levels of in-flight disturbance were considered as reported in Table~\ref{tab:inflight_noise_levels} and called Low (L), Medium (M), and High (H).
In each case, the standard deviation of the noise in the thrust direction is 5 times greater than the one in a random direction, thus $\alpha_{G_\parallel} = 5$. 

\begin{table}[h]
    \caption{Standard deviations of the random-direction Gaussian in-flight disturbance in terms of $T/m$ percentage}
    \label{tab:inflight_noise_levels}
    \centering
    \begin{tabular}{c c c c}
        \toprule
        \bf Parameter & \bf Case L & \bf Case M & \bf Case H \\ 
        \midrule
        $\alpha_G$ & \SI{0.25}{\permille}  & \SI{0.5}{\permille} & \SI{1}{\permille} \\
        \bottomrule
    \end{tabular}
\end{table}

In each simulation, the additional burn time $\Delta t_b$ is sampled uniformly in the range $[0, 5]$\si{\s}, while the additional total impulse $\Delta I_{\text{tot}}$ is a fraction of the nominal one uniformly sampled in the range $[0, 5]$\si{\permille}.
The update frequency of the MPC is set to \SI{1}{\Hz}, meaning that the OCP is solved every $T = \SI{1}{\s}$.

\begin{table*}
    \begin{center}
        \caption{Results of the Monte Carlo campaigns}
        \label{tab:MC_results}
        \begin{tabular}{c c c c c c c c c c}
        \toprule
        \multirow{2}[2]{*}{\shortstack[c]{\bf Algorithm}} &
        \multirow{2}[2]{*}{\shortstack[c]{\bf Case}} &
        \multicolumn{4}{c}{\bf Payload mass (kg)}   &
        \multicolumn{4}{c}{\bf $\varphi_{R}$ (deg)}   \\
        \cmidrule(lr){3-6} 
        \cmidrule(lr){7-10}
        & & Min & Mean & Max & $\sigma$ & Min & Mean & Max & $\sigma$ \\
        \midrule
        \multirow{3}[0]{*}{\shortstack[c]{Single \\ return}}
        & L &
        \num{1359.19} & \num{1409.74} & \num{1449.70} & \num{17.10} & \num{59.74} & \num{63.60} & \num{68.15} & \num{1.87} \\
        & M &
        \num{1335.17} & \num{1405.21} & \num{1482.21} & \num{24.41} & \num{59.30} & \num{63.33} & \num{68.65} & \num{2.00} \\
        & H &
        \num{1254.00} & \num{1407.49} & \num{1503.37} & \num{39.77} & \num{57.78} & \num{63.63} & \num{75.40} & \num{2.73} \\
        \midrule
        \multirow{3}[0]{*}{\shortstack[c]{Multi \\ return}}
        & L & 
        \num{1278.39} & \num{1331.44} & \num{1375.78} & \num{18.16} & \num{59.72} & \num{60.00} & \num{60.21} & \num{0.09} \\
        & M & 
        \num{1266.97} & \num{1329.81} & \num{1389.62} & \num{23.56} & \num{59.17} & \num{60.01} & \num{60.62} & \num{0.20} \\
        & H & 
        \num{1215.90} & \num{1331.64} & \num{1456.62} & \num{41.06} & \num{58.25} & \num{60.02} & \num{61.41} & \num{0.38} \\
        \bottomrule
        \end{tabular} 
    \end{center}
\end{table*}

The results of the Monte Carlo campaigns are summarized in Table~\ref{tab:MC_results}.
For each case, 400 independent MPC simulations were carried out.
For each simulation, the final conditions at the burnout of the SRM are used to estimate the attainable payload mass by solving the ascent OCP from Phases \CoastingThree{} to \StageFourTwo{} and to evaluate the splash-down point through forward propagation of the equations of motion, Eqs.~\eqref{eq:original_ODE_r}--\eqref{eq:original_ODE_m}.
As expected, for both algorithms, as the intensity of the in-flight disturbance increases, larger dispersions on the payload mass and on the splash-down point are observed.
The average payload mass attained by the algorithms is slightly greater than the nominal one due to the additional $I_{\text{tot}}$ that provides more energy to the system.
The gain in payload mass of the multi-return strategy is smaller due to the fact that the additional impulse is provided in the neutral axis direction, which is not optimal.
Instead, the payload mass ranges of the two algorithms are comparable, except for a shift of approximately \SI{75}{\kilo\g} due to the cost of the neutral axis maneuver.

The major difference in the behavior of the two strategies lies in the splash-down point.
The additional impulse provided by the SRM shifts the mean splash-down point attained by the single-return strategy farther by approximately \SI{3.6}{\degree} in all three cases, meaning that this deviation depends only on the SRM performance uncertainty.
Instead, the multi-return algorithm retains the mean splash-down location with great accuracy (the largest error is \SI{0.02}{\degree} in case H) thanks to the fact that the thrust vector is oriented in the neutral axis direction during the extra seconds of operation.
The dispersion of the spent stage impact latitude grows for both approaches as the disturbance intensity increases, but, while the lower bounds are comparable, the upper bounds of the return points corresponding to the single-return strategy are much more distant from the desired value.
This asymmetry is due to the fact that the additional SRM impulse moves the return point only farther; thus, the two disturbances add up for the single-return approach.
Instead, the neutral axis maneuver compensates for the SRM performance uncertainty and the range of the attained $\varphi_R$ is symmetric.
Also, the dispersion of the return point is significantly reduced when employing the multi-return algorithm, with standard deviations $\sigma$ one order of magnitude smaller than in the single-return approach.



\begin{figure}
    \centering

    \begin{minipage}{\columnwidth}
        \centering
        \includegraphics[width=1\linewidth]{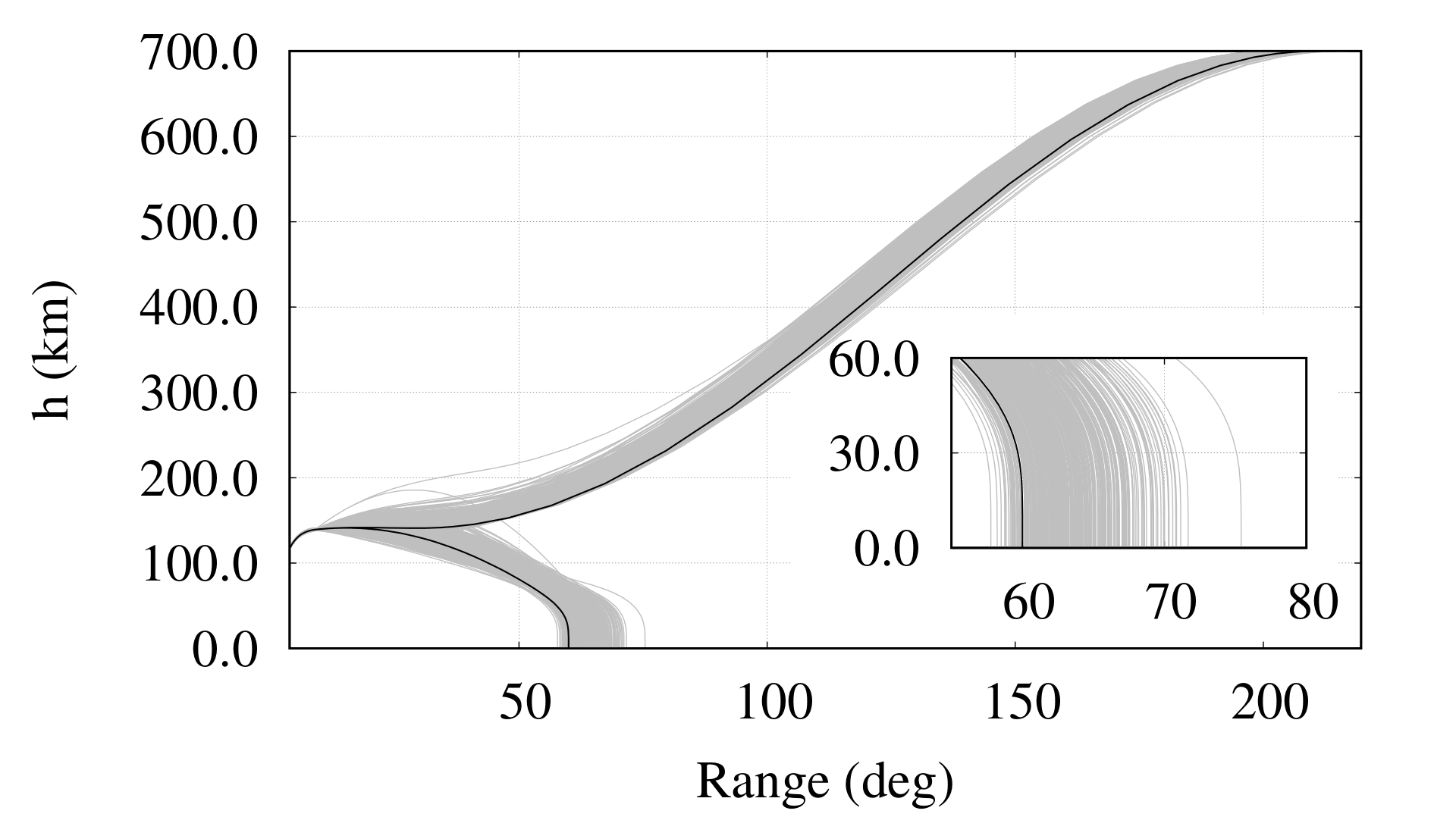}
        \subcaption{Single-return algorithm}
        \label{fig:MC_altitudes_single}
    \end{minipage}
    
\centering
    \begin{minipage}{\columnwidth}
        \centering
        \includegraphics[width=1\linewidth]{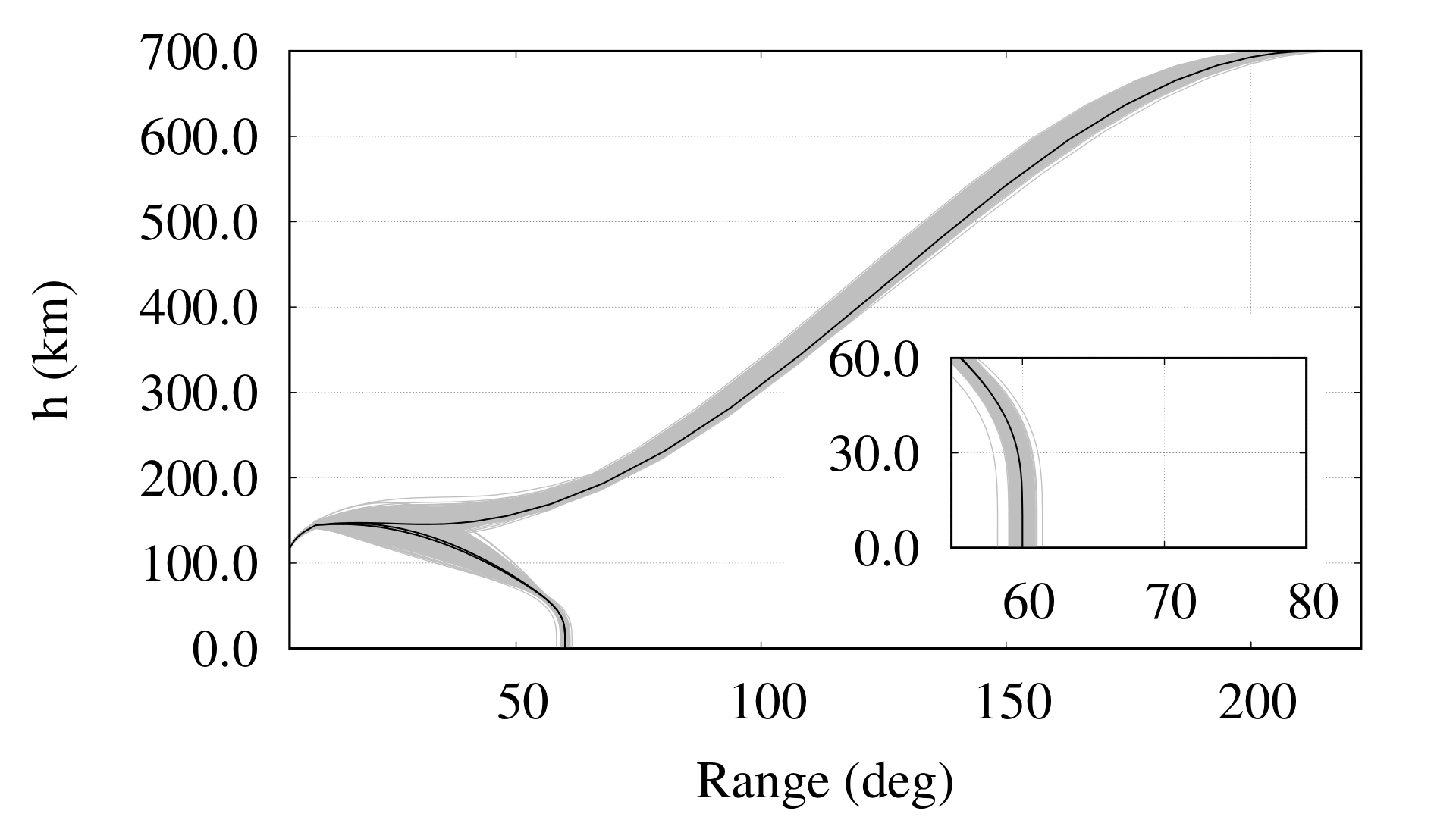}
        \subcaption{Multi-return algorithm}
        \label{fig:MC_altitudes_multi}
    \end{minipage}
    
    \caption{Altitude profiles for Case H}
    \label{fig:MC_altitudes}
\end{figure}

The trajectories for case H, which is associated with the largest envelopes, are reported in Fig.~\ref{fig:MC_altitudes} and compared with the nominal one, which is denoted by the black line.
The trajectories of both guidance strategies deviate significantly from the nominal trajectory during the central portion of the flight, as the MPC approach autonomously recomputes optimal paths to efficiently compensate for the encountered disturbances.
However, while the terminal conditions at the payload release are met by both algorithms, the single-return strategy features a much greater dispersion on the return point, as previously commented.

\begin{figure}[h]
    \centering

    \begin{minipage}{\columnwidth}
        \centering
        \includegraphics[width=1\linewidth]{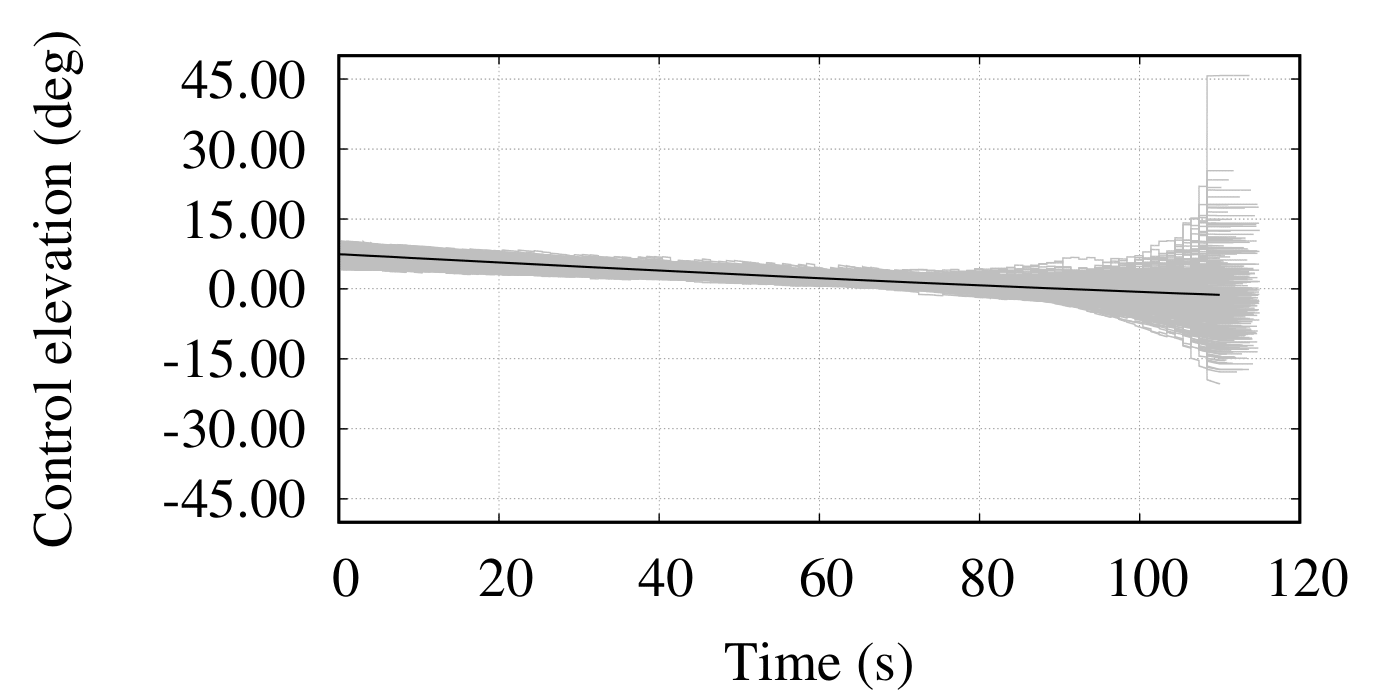}
        \subcaption{Single-return algorithm}
        \label{fig:MC_elevations_single}
    \end{minipage}
    
\centering
    \begin{minipage}{\columnwidth}
        \centering
        \includegraphics[width=1\linewidth]{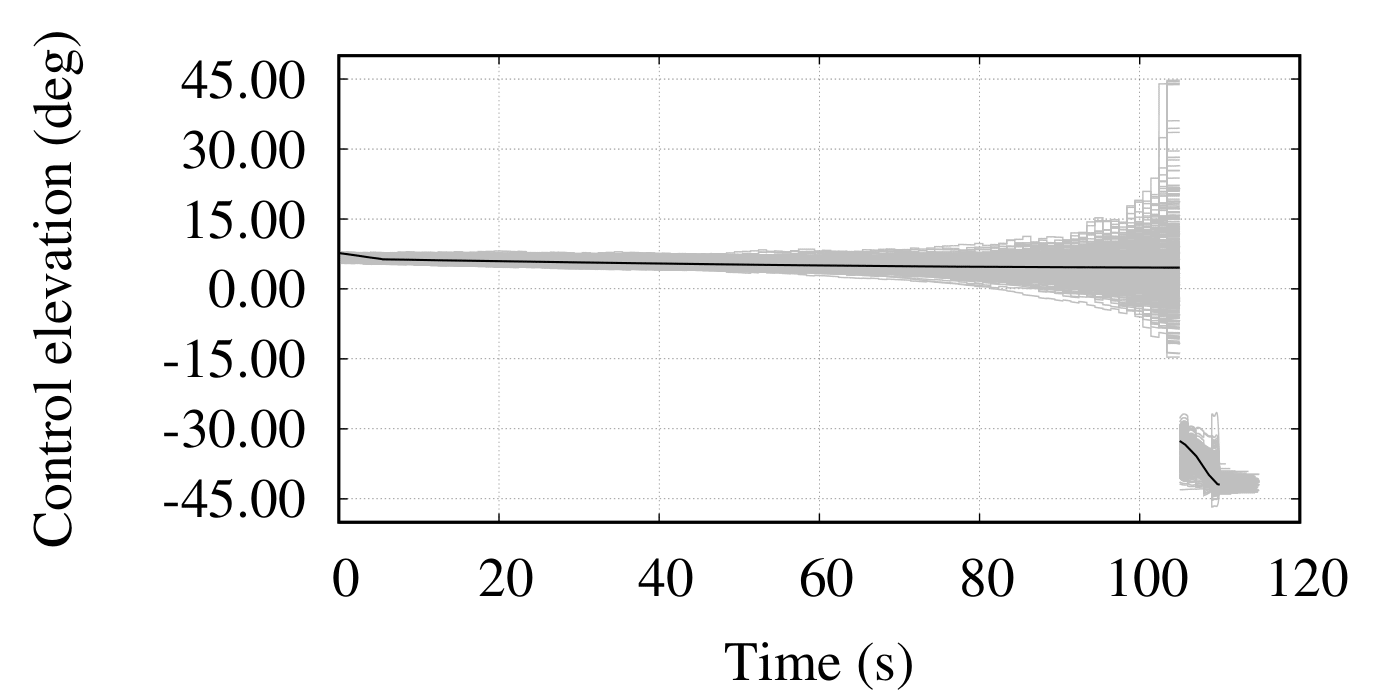}
        \subcaption{Multi-return algorithm}
        \label{fig:MC_elevations_multi}
    \end{minipage}
    
    \caption{Control elevation profiles for Case H}
    \label{fig:MC_elevations}
\end{figure}

The control histories are shown in Fig.~\ref{fig:MC_elevations} in terms of the elevation angle, which is defined as the angle between the thrust direction $\bm{\hat{T}}$ and the local horizontal.
The elevation profiles are close to the nominal controls for most of the third stage operation, with larger deviations toward the end.
Indeed, due to the external disturbance, meeting the splash-down constraint requires rapid variations of the thrust direction as the burnout time approaches.
Figure~\ref{fig:MC_elevations_multi} shows the different orientation of the neutral axis direction compared to the optimal one, as the maneuver requires the rocket to quickly rotate by approximately \SI{40}{\degree} downward.

\begin{figure}[h]
    \centering

    \begin{minipage}{0.9\columnwidth}
        \centering
        \includegraphics[width=1\linewidth]{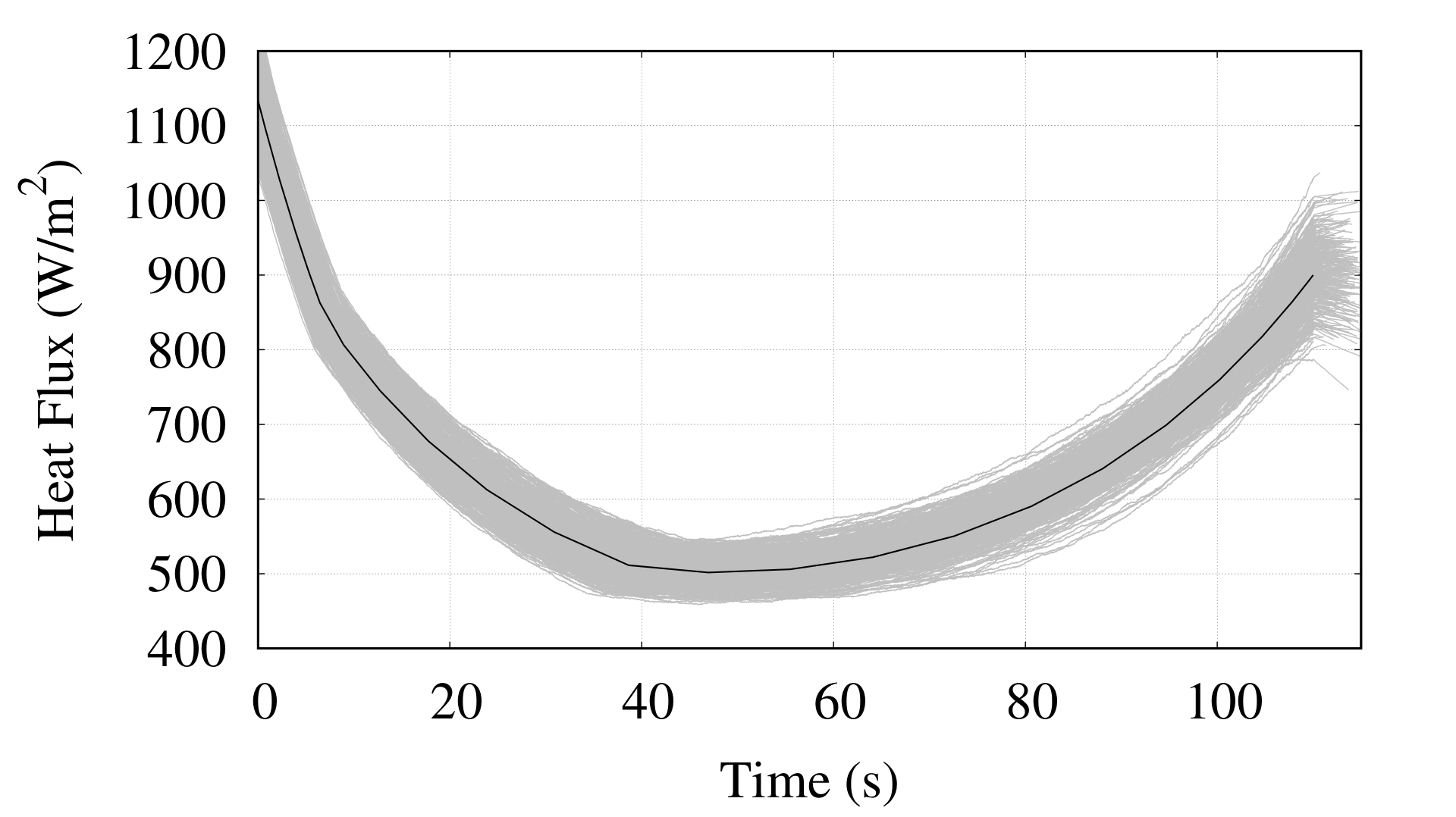}
        \subcaption{Single-return algorithm}
        \label{fig:MC_hfs_single}
    \end{minipage}
    
    \centering
    \begin{minipage}{0.9\columnwidth}
        \centering
        \includegraphics[width=1\linewidth]{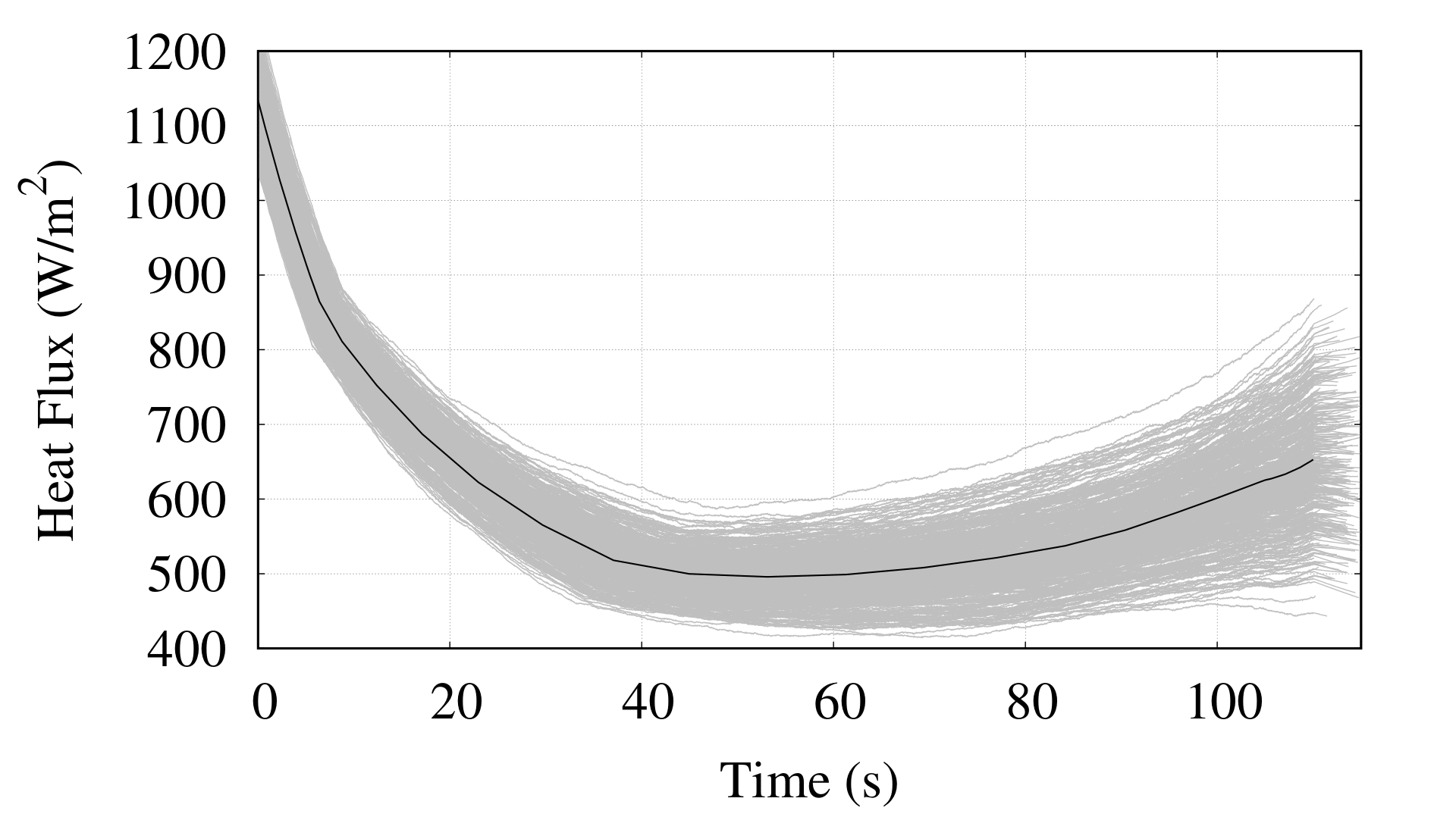}
        \subcaption{Multi-return algorithm}
        \label{fig:MC_hfs_multi}
    \end{minipage}
    
    \caption{Heat flux profiles for Case H}
    \label{fig:MC_hfs}
\end{figure}

The heat flux profiles are illustrated in Fig.~\ref{fig:MC_hfs}.
Despite the relaxation, the thermal requirement is satisfied within an acceptable error.
In particular, the heat flux is never greater than \SI{1000}{\watt\per\square\m}, except for a few simulations of the single-return strategy in the proximity of the burnout where, due to the additional seconds of SRM operations, the relative velocity increases more than predicted.
However, the violation only depends on the scattering of the initial conditions, as increased initial velocity and reduced altitude expose the system to more critical thermal conditions. 


As a final remark, 
the average solution time of the single-return OCP was \SI{1.50}{\s}, while solving the multi-return problem required \SI{3.20}{\s} on average.
The multi-return OCP is thus more expensive to solve, due to the additional return phase and the related constraints that significantly increase the problem dimension and complexity.
Even though both times are incompatible with a \SI{1}{\Hz} update frequency, this issue could be easily addressed by using optimized software, dedicated hardware, and custom convex solvers, for which the computation times are expected to be at least one order of magnitude lower.
Also, previous investigations on this topic show that good MPC performance could be obtained even with less frequent updates, in the range \num{1}--\SI{3}{\s} \cite{benedikter2020autonomous}.

\section{Conclusions}
\label{sec:conclusions}

This paper presented a model predictive control algorithm for the optimal guidance of a launch vehicle upper stage that robustly controls the splash-down location of the spent stage even in presence of uncertainties on the motor burn time, by incorporating in the OCP a set of boundary conditions that allows to numerically find the thrust direction corresponding to the null miss condition.
Specifically, two return phases of the spent stage are considered in every OCP, a nominal return and another perturbed by a small velocity impulse that is applied in the same direction as the thrust at the end of the nominal cut-off.
This approach avoids the solution of the neutral axis equation, which enforces the same condition algebraically for a ballistic return phase, and can be directly introduced in the MPC guidance scheme, 
thus greatly reducing the validation and verification burden  of such a complex system.

The OCP was mindfully cast as a sequence of convex problems that converges to the optimal solution in a short time, allowing the MPC algorithm to attain an high update frequency and effectively compensate for external random disturbances.
Also, a $hp$-pseudospectral discretization was adopted to provide an accurate representation of the continuous-time ascent and return dynamics with a limited computational burden.

A Monte Carlo analysis has been carried out to investigate the effectiveness of the proposed approach in presence of in-flight disturbances, off-nominal operating conditions, and uncertain SRM performance.
The results of this analysis provide significant evidence of the robustness of the MPC algorithm.
In particular, the benefits in terms of robustness of the splash-down location of the present approach compared to the previous one that did not incorporate a null miss constraint into the OCP are apparent.
Indeed, even though the magnitude of the in-flight disturbances directly affects the dispersion on the payload mass and on the splash-down location, the deviations are one order of magnitude smaller.

Even though further tests are still necessary to demonstrate the suitability of the proposed MPC algorithm for practical onboard implementation, the preliminary results are very encouraging, as the design of the neutral axis maneuver is still a burdensome, yet necessary, task for several launch vehicles, first and foremost, VEGA, and having an autonomous strategy to include the null miss condition within the optimization loop represents a significant step forward in the design of launch vehicles guidance and control algorithms. 

\bibliographystyle{ieeetr}
\bibliography{references}

\end{document}